\documentclass[12pt]{article}
\usepackage{amsthm,amsfonts,amssymb,amscd}

\begin{document}

\begin{center}
\Large\bf On the $8n^2$-inequality
\end{center}
\vspace{0.5cm}

\begin{center}
\bf Aleksandr Pukhlikov
\end{center}
\vspace{0.5cm}

\parshape=1
3cm 10cm \noindent {\small \quad \quad \quad
\quad\quad\quad\quad\quad\quad\quad {\bf Abstract}\newline We give
a complete proof of the so called $8n^2$-inequality, a local
inequality for the self-intersection of a movable linear system at
an isolated centre of a non canonical singularity. The inequality
was suggested and several times published by I.Cheltsov but some
of his arguments are faulty. We explain the mistake and replace
the faulty piece by a correct argument.} \vspace{1cm}

{\bf 1. Introduction.} The aim of this note is to give a complete
proof of the so called $8n^2$-inequality, correcting the mistakes
in [1,2] and some other papers where the erroneous arguments were
reproduced. That inequality makes it possible to prove birational
(super)rigidity of several types of Fano varieties of
anticanonical degrees 6,7 and 8. One step in the arguments of
[1,2] is based on an erroneous claim which it is unclear how to
correct (and whether that step can be corrected at all, following
the approach of [1,2]).

In this note, we replace the faulty step by a different argument,
thus making the proof complete.

The note is organized as follows. In Sec.2 we formulate the
$8n^2$-inequality and reproduce briefly that part of the proof in
[1,2], which is correct. In Sec. 3 we present the new arguments,
completing the proof. In Sec. 4 we discuss the mistakes in [1,2],
one of which is really serious and undermines the whole proof. The
list of references given in Sec. 5 is far from being complete: the
faulty arguments were published in a few more papers.

It is worth mentioning that in the survey [3] the arguments are
the same as those of [1,2], with one exception: the very claim,
for which [1,2] give the faulty arguments, is presented in [3]
with essentially no proof at all.\vspace{0.3cm}

{\bf 2. The inequality and start of the proof.} Let $o\in X$ be a
germ of a smooth variety of dimension $\mathop{\rm dim} X\geq 4$.
Let $\Sigma$ be a movable system on $X$ and the effective cycle
$$
Z=(D_1\circ D_2),
$$
where $D_1,D_2\in Z$ are generic divisors, its self intersection.
Let us blow up the point $o$ on $X$:
$$
\varphi\colon X^+\to X,
$$
$E=\varphi^{-1}(o)\cong{\mathbb P}^{\rm dimX-1}$ is the
exceptional divisor. The strict transforms of the system $\Sigma$
and the cycle $Z$ on $X^+$ denote by the symbols $\Sigma^+$ and
$Z^+$, respectively.

{\bf Theorem ($8n^2$-inequality)}. {\it Assume that the pair
$$
(X,\frac{1}{n}\Sigma)
$$
is not canonical, but canonical outside the point $o$, where $n$
is some positive integer. There exists a linear subspace $P\subset
E$ of codimension two (with respect to $E$) such, that the
following inequality holds:}
$$
\mathop{\rm mult}\nolimits_oZ+\mathop{\rm
mult}\nolimits_PZ^+>8n^2.
$$

{\bf Proof.} Note that if $\mathop{\rm mult}_oZ>8n^2$, then for
$P$ we can take any subspace of codimension two in $E$.

We start the proof arguing as in [1,2].

Restricting $\Sigma$ onto a germ of a generic smooth subvariety,
containing the point $o$, we may assume that $\mathop{\rm
dim}X=4$. Moreover, we may assume that $\nu=\mathop{\rm
mult}_oZ\leq 2\sqrt{2}n< 3n$, since otherwise
$$
\mathop{\rm mult}\nolimits_oZ\geq\nu^2>8n^2
$$
and there is nothing to prove.

{\bf Lemma 1.} {\it The pair
\begin{equation}\label{x1}
(X^+,\frac{1}{n}\Sigma^++\frac{(\nu-2n)}{n}E)
\end{equation}
is not log canonical, and the centre of any of its non log
canonical singularities is contained in the exceptional divisor
$E$.}

{\bf Proof.} Let $\lambda\colon\widetilde{X} \to X$ be a
resolution of singularities of the pair $(X,\frac{1}{n}\Sigma)$
and $E^*\subset \widetilde{X}$ a prime exceptional divisor,
realizing a non-canonical singularity of that pair. Then
$\lambda(E^*)=o$ and the Noether-Fano inequality holds:
$$
\nu_{E^*}(\Sigma)>na(E^*).
$$
For a generic divisor $D\in\Sigma$ we get $\varphi^*D=D+\nu E$, so
that
$$
\nu_{E^*}(\Sigma)=\nu_{E^*}(\Sigma^+)+\nu\cdot\nu_{E^*}(E)
$$
and
$$
a(E^*,X)=a(E^*,X^+)+3\nu_{E^*}(E).
$$
From here we get
$$
\nu_{E^*}\left(\frac{1}{n}\Sigma^++\frac{\nu-2n}{n}E\right)=
\nu_{E^*}\left(\frac{1}{n}\Sigma\right)-2\nu_{E^*}(E)>
$$
$$
>a(E^*,X^+)+\nu_{E^*}(E)\geq a(E^*,X^+)+1,
$$
which proves the lemma.

Let $R\ni o$ be a generic three-dimensional germ, $R^+\subset X^+$
its strict transform on the blow up of the point $o$. For a small
$\varepsilon>0$ the pair
$$
\left(X^+,\frac{1}{1+\varepsilon}\frac{1}{n}\Sigma^+
+\frac{\nu-2n}{n}E+R^+\right)
$$
still satisfies the connectedness principle (with respect to the
morphism $\varphi\colon X^+\to X$), so that the set of centres of
non log canonical singularities of that pair is connected. Since
$R^+$ is a non log canonical singularity itself, we obtain, that
there is a non log canonical singularity of the pair (\ref{x1}),
the centre of which on $X^+$ is of positive dimension, since it
intersects $R^+$.

Let $Y\subset E$ be a centre of a non log canonical singularity of
the pair (\ref{x1}) that has the maximal dimension.

If $\mathop{\rm dim}Y=2$, then consider a generic two-dimensional
germ $S$, intersecting $Y$ transversally at a point of general
position. The restriction of the pair (\ref{x1}) onto $S$ is not
log canonical at that point, so that, arguing as in [1,2], we see
that
$$
\mathop{\rm mult}\nolimits_Y(D^+_1\circ
D^+_2)>4\left(3-\frac{\nu}{n}\right)n^2,
$$
so that
$$
\mathop{\rm mult}\nolimits_oZ\geq\nu^2+\mathop{\rm
mult}\nolimits_Y(D^+_1\circ D^+_2)\mathop{\rm deg}Y>
$$
$$
>(\nu-2n)^2+8n^2,
$$
which is what we need.

If $\mathop{\rm dim}Y=1$, then, since the pair
\begin{equation}\label{x2}
\left(R^+,\frac{1}{1+\varepsilon}\frac{1}{n}\Sigma^+_R
+\frac{\nu-2n}{n}E_R\right),
\end{equation}
where $\Sigma^+_R=\Sigma^+|_{R^+}$ and $E_R=E|_{R^+}$, satisfies
the condition of the connectedness principle and $R^+$ intersects
$Y$ at $\mathop{\rm deg}Y$ distinct points, we conclude that
$Y\subset E$ is a line in ${\mathbb P}^3$.

Now we need to distinguish between the following two cases: when
$\nu\geq 2n$ and when $\nu<2n$. The methods of proving the
$8n^2$-inequality in these two cases are absolutely different.
Consider first the case $\nu\geq 2n$.

Let us choose as $R\ni o$ a generic three-dimensional germ,
satisfying the condition $R^+\supset Y$. Since the pair (\ref{x2})
is effective (recall that $\nu\geq 2n$), one may apply inversion
of adjunction [4, Chapter 17] and conclude that the pair
(\ref{x2}) is not log canonical at $Y$.

Now arguing in the same way as for $\mathop{\rm dim}Y=2$, with
$R^+\supset Y$, we get the inequality
$$
\mathop{\rm mult}\nolimits_Y(D^+_1|_{R^+}\circ
D^+_2|_{R^+})>4\left(3-\frac{\nu}{n}\right)n^2.
$$
On the left in brackets we have the self-intersection of the
movable system $\Sigma^+_R$, which breaks into two natural
components:
$$
(D^+_1|_{R^+} \circ D^+_2|_{R^+})=Z^+_R+Z^{(1)}_R,
$$
where $Z^+_R$ is the strict transform of the cycle $Z_R=Z|_R$ on
$R^+$ and the support of the cycle $Z^{(1)}_R$ is contained in
$E_R$. The line $Y$ is a component of the effective cycle
$Z^{(1)}_R$.

On the other hand, for the self-intersection of the movable linear
system $\Sigma^+$ we get
$$
(D^+_1\circ D^+_2)=Z^++Z_1,
$$
where the support of the cycle $Z_1$ is contained in $E$. From the
genericity of $R$ it follows that outside the line $Y$ the cycles
$Z^{(1)}_R$ and $Z_1|_{R^+}$ coincide, whereas for $Y$ we get the
equality
$$
\mathop{\rm mult}\nolimits_YZ^{(1)}_R=\mathop{\rm
mult}\nolimits_YZ^++\mathop{\rm mult}\nolimits_YZ_1.
$$
However, $\mathop{\rm mult}_YZ_1\leq\mathop{\rm deg}Z_1$, so that
$$
\mathop{\rm mult}\nolimits_oZ+\mathop{\rm mult}\nolimits_YZ^+=
$$
$$
=\nu^2+\mathop{\rm deg}Z_1+\mathop{\rm mult}\nolimits_YZ^+\geq
$$
$$
\geq\nu^2+\mathop{\rm mult}\nolimits_YZ^{(1)}_R>8n^2,
$$
which is what we need. This completes the case $\nu\geq 2n$.

Note that the key point in this argument is that the pair
(\ref{x2}) is effective. For $\nu<2n$  inversion of adjunction can
not be applied (as it was done in [3]). The additional arguments
in [1,2], proving inversion of adjunction specially for this pair
for $\nu<2n$, are faulty.\vspace{0.3cm}

{\bf 3. The technique of counting multiplicities.} Starting from
this moment, assume that $\nu<2n$.

Consider again the pair (\ref{x2}) for a generic germ $R\ni o$.
Let $y=Y\cap R^+$ be the point of (transversal) intersection of
the line $Y$ and the variety $R^+$. Since $a(E_R,R)=2$, the non
log canonicity of the pair (\ref{x2}) at the point $y$ implies the
non log canonicity of the pair
$$
\left(R,\frac{1}{n}\Sigma_R\right)
$$
at the point $o$, whereas the centre of some non log canonical
(that is, log maximal) singularity on $R^+$ is a point $y$.

Now the $8n^2$-inequality comes from the following fact.

{\bf Lemma 2.} {\it The following inequality holds:
$$
\mathop{\rm mult}\nolimits_oZ_R+\mathop{\rm
mult}\nolimits_yZ^+_R>8n^2,
$$
where $Z_R$ is the self-intersection of a movable linear system
$\Sigma_R$ and $Z^+_R$ is its strict transform on $R^+$.}

{\bf Proof.} Consider the resolution of the maximal singularity of
the system $\Sigma_R$, the centre of which on $R^+$ is the point
$y$:
$$
\begin{array}{ccc}
R_i & \stackrel{\psi_i}{\to} & R_{i-1}\\
\cup & & \cup\\
E_i & & B_{i-1},\\
\end{array}
$$
where $B_{i-1}$ is the centre of the singularity on $R_{i-1}$,
$E_i=\psi^{-1}_i(B_{i-1})$ is the exceptional divisor, $B_0=o$,
$B_1=y\in E_1$, $i=1,\dots,N$, where the first $L$ blow ups
correspond to points, for $i\geq L+1$ curves are blown up. Since
$$
\mathop{\rm mult}\nolimits_o\Sigma_R=\mathop{\rm
mult}\nolimits_o\Sigma<2n,
$$
we get $L<N$, $B_L\subset E_L\cong{\mathbb P}^2$ is a line and for
$i\geq L+1$
$$
\mathop{\rm deg}[\psi_i|_{B_i}\colon B_i\to B_{i-1}]=1,
$$
that is, $B_i\subset E_i$ is a section of the ruled surface $E_i$.

Consider the graph of the sequence of blow ups $\psi_i$.

{\bf Lemma 3.} {\it The vertices $L+1$ and $L-1$ are not connected
by an arrow:}
$$
L+1\nrightarrow L-1.
$$

{\bf Proof.} Assume the converse: $L+1\to L-1$. This means that
$$
B_L=E_L\cap E^L_{L-1}
$$
is the exceptional line on the surface $E^L_{L-1}$ and the map
$$
E^{L+1}_{L-1}\to E^L_{L-1}
$$
is an isomorphism. As usual, set
$$
\nu_i=\mathop{\rm mult}\nolimits_{B_{i-1}}\Sigma^{i-1}_R,
$$
$i=1,\dots,N$. Let us restrict the movable linear system
$\Sigma^{L+1}_R$ onto the surface $E^{L+1}_{L-1}$ (that is, onto
the plane $E_{L-1}\cong{\mathbb P}^2$ with the blown up point
$B_{L-1}$). We obtain a non-empty (but, of course, not necessarily
movable) linear system, which is a subsystem of the complete
linear system
$$
\left|\nu_{L-1}(-E_{L-1}|_{E_{L-1}})-(\nu_L+\nu_{L+1})B_L\right|.
$$
Since $(-E_{L-1}|_{E_{L-1}})$ is the class of a line on the plane
$E_{L-1}$, this implies that
$$
\nu_{L-1}\geq\nu_L+\nu_{L+1}>2n,
$$
so that the more so $\nu_1=\nu>2n$. A contradiction. Q.E.D. for
the lemma.

Set, as usual,
$$
m_i=\mathop{\rm mult}\nolimits_{B_{i-1}}(Z_R)^{i-1},
$$
$i=1,\dots,L$, so that, in particular,
$$
m_1=\mathop{\rm mult}\nolimits_oZ_R\quad\mbox{and}\quad
m_2=\mathop{\rm mult}\nolimits_yZ^+_R.
$$
Let $p_i\geq 1$ be the number of paths in the graph of the
sequence of blow ups $\psi_i$ from the vertex $N$ to the vertex
$i$, and $p_N=1$ by definition, see [5,6]. By what we proved,
$$
p_N=p_{N-1}=\dots=p_L=p_{L-1}=1,
$$
and the number of paths $p_i$ for $i\leq L$ is the number of paths
from the vertex $L$ to the vertex $i$. By the technique of
counting multiplicities [5,6], we get the inequality
$$
\sum^L_{i=1}p_im_i\geq\sum^N_{i=1}p_i\nu^2_i
$$
and, besides, the Noether-Fano inequality holds:
$$
\sum^N_{i=1}p_i\nu_i>
n\left(2\sum^L_{i=1}p_i+\sum^N_{i=L+1}p_i\right).
$$
(In fact, a somewhat stronger inequality holds, the {\it log}
Noether-Fano inequality, but we do not need that.) From the last
two estimates one obtains in the standard way [5,6] the inequality
$$
\sum^L_{i=1}p_im_i>\frac{(2\Sigma_0+\Sigma_1)^2}
{\Sigma_0+\Sigma_1}n^2,
$$
where $\Sigma_0=\sum\limits^L\limits_{i=1}p_i$ and
$\Sigma_1=\sum\limits^N\limits_{i=L+1}p_i=N-L$. Taking into
account that for $i\geq 2$ we get
$$
m_i\leq m_2
$$
and the obvious inequality
$(2\Sigma_0+\Sigma_1)^2>4\Sigma_0(\Sigma_0+\Sigma_1)$, we obtain
the following estimate
$$
p_1m_1+(\Sigma_0-p_1)m_2>4n^2\Sigma_0.
$$
Now assume that the claim of the lemma is false:
$$
m_1+m_2\leq 8n^2.
$$

{\bf Lemma 4.} {\it The following inequality holds:} $\Sigma_0\geq
2p_1$.

{\bf Proof.} By definition,
$$
p_1=\sum_{i\to 1}p_i,
$$
however, by Lemma 3 from $i\to 1$ it follows that $i\leq L$, so
that $p_1\leq\Sigma_0-p_1$, which is what we need. Q.E.D. for the
lemma.

Now, taking into account that $m_2\leq m_1$, we obtain
$$
p_1m_1+(\Sigma_0-p_1)m_2=p_1(m_1+m_2)+(\Sigma_0-2p_1)m_2\leq
$$
$$
\leq 8p_1n^2+(\Sigma_0-2p_1)\cdot 4n^2=4n^2\Sigma_0.
$$
This is a contradiction. Q.E.D. for Lemma 2.

Proof of our theorem is complete.\vspace{0.3cm}

{\bf 4. On the faulty arguments.} The proof of the
$8n^2$-inequality given in [1] is invalid. We refer the reader to
that paper (the numbers of pages and claims correspond to the
archive version given in the reference [1]).\vspace{0.1cm}

{\bf Mistake 1.} Page 6 in the archive version of [1], just after
the proof of Lemma 27. The claim that the intersection of the
divisor $S$ with the curve $C$ is either trivial or consists of
more than one point, is wrong. The divisor $F$ contains a
3-dimensional family of smooth rational curves, intersecting
transversally a generic divisor дивизор $S$ at one point. Namely,
the surface
$$
\bar{E}\cap F
$$
(in $E\cong{\mathbb P}^3$ we blow up the line $L$, and
$\bar{E}\cap F$ is the exceptional divisor) is ${\mathbb
P}^1\times{\mathbb P}^1$, so that any curve $C$ of bidegree (1,1)
can be used. (This error can be corrected, for instance, by
restricting the linear system onto the exceptional divisor $E$ and
showing that in the case of a $(1,1)$ curve $\nu>2n$, so that the
arguments based on inversion of adjunction work.) \vspace{0.1cm}

{\bf Mistake 2.} Page 5 of [1]: Corollary 24 is not true. The
intersection $S\cap C$ can well be empty.

Indeed, from the following two facts:\vspace{0.1cm}

{\bf A.} the set $LCS(S,(B^W+\bar{E}+2F)|_S)$ either consists of
one point or contains a curve (which is deduced from the
connectedness principle (Theorem 14)), and\vspace{0.1cm}

{\bf B.} there is a curve $C$, a section of the bundle $F\to L$,
such that $C$ is the unique element of the set
$LCS(W,B^W+\bar{E}+aF)$, $a=1,2$, which is contained in $F$ and
dominates $L$ (just above, page 5),\vspace{0.1cm}

\noindent it {\bf does not} follow that
$LCS(S,(B^W+\bar{E}+2F)|_S)$ is the intersection $S\cap C$ (and
this is exactly how Corollary 24 is proved), since
$LCS(S,(B^W+\bar{E}+2F|_S))$ can well be the intersection of $S$
with the centre of a non log canonical singularity of the pair
$(W,B^W+\bar{E}+2F)$, which {\bf does not} dominate $L$: for
instance, with a line in the fiber of the bundle $F\to L$ (this is
a ${\mathbb P}^2$-bundle over $L$). In particular, if the set
$LCS(W,B^W+\bar{E}+aF)$, $a=1,2$, is a connected union of two
curves:\vspace{0.1cm}

(1) a line in a fiber $F\to L$ and\vspace{0.1cm}

(2) a curve of bidegree (1,0) on the surface $\bar{E}\cap F$,
which is ${\mathbb P}^1\times{\mathbb P}^1$, that is, a section of
$F\to L$, contained in $\bar{E}\cap F$ and having the zero
self-intersection on that surface,\vspace{0.1cm}

\noindent then A and B hold, but there is no contradiction at all.

The ``proof'' of Corollary 24 is faulty, because the (correct)
claim ``the centre of any singularity, {\it dominating} $L$, is
$C$'' is used actually as the claim ``the centre of any
singularity is $C$''. The example above shows that the arguments
of [1] are faulty and give no proof of the $8n^2$-inequality. (In
[2] the same arguments are given as those used in [1] for proving
B, after which it is claimed that $C$ is the unique element of the
set $LCS(\dots)$, without mentioning that $L$ is dominated. Here
it is easier to see the point of trouble.) \vspace{0.3cm}

\newpage

{\bf 5. References.}

\noindent [1] Cheltsov I., Double cubics and double quartics,
Math. Z. 253 (2006), no. 1, 75-86; arXiv:math/0410408v4 [math.AG]
\vspace{0.1cm}

\noindent [2] Cheltsov I., Non-rationality of a four-dimensional
smooth complete intersection of a quadric and a quartic, not
containing a plane, Sbornik: Mathematics, 194 (2003),
1679-1699.\vspace{0.1cm}

\noindent[3] Iskovskikh V.A., Birational rigidity of Fano
hypersurfaces in the framework of Mori theory, Russian Math.
Surveys, 2001, V. 56, No. 2.\vspace{0.1cm}

\noindent [4] Koll{\'a}r J., et al., Flips and Abundance for
Algebraic Threefolds, Asterisque 211, 1993. \vspace{0.1cm}

\noindent [5] Pukhlikov A.V., Essentials of the method of maximal
singularities, in ``Explicit Birational Geometry of Threefolds'',
London Mathematical Society Lecture Note Series {\bf 281} (2000),
Cambridge University Press, 73-100. \vspace{0.1cm}

\noindent [6] Pukhlikov A.V., Birationally rigid varieties. I.
Fano varieties. Russian Math. Surveys. 2007. V. 62, no. 5,
857-942. \vspace{0.1cm}

\end{document}